\documentclass{elsart}
\usepackage{color}

\usepackage{amssymb}
\usepackage{amsmath}
\usepackage{latexsym}
\usepackage{amsfonts}
\usepackage{amssymb}
\usepackage{rotating}
\usepackage{soul}

\newcommand{\essi}{\operatornamewithlimits{ess\,inf}}
\newcommand{\esss}{\operatornamewithlimits{ess\,sup}}

\DeclareMathOperator{\rota}{rot}

\newcommand{\al}{\alpha}

\newcommand{\dl}{\delta}

\newcommand{\Om}{\Omega}
\newcommand{\lb}{\lambda}

\newcommand{\ve}{\varepsilon}
\newcommand{\gm}{\gamma}

\newcommand{\vi}{\varphi}

\newcommand{\intl}{\int\limits}

\newcommand{\ol}{\overline}

\begin{document}

\begin{frontmatter}



\title{HARDY TYPE INEQUALITY IN \\ VARIABLE LEBESGUE SPACES\thanksref{label1}}
  \thanks[label1]{This work was made under the project \textit{Variable Exponent Analysis} supported by INTAS grant
 Nr.06-1000017-8792.}
 \thanks[label2]{ Supported by \textit{Funda\c c\~ao para a
Ci\^encia e a Tecnologia} (FCT) (Grant Nr. SFRH/BD/22977/2005) of the Portuguese Government.}
\author[Algarve]{Humberto Rafeiro\thanksref{label2}},
\ead{hrafeiro@ualg.pt}
\ead[url]{w3.ualg.pt/\textasciitilde  hrafeiro}
\author[Algarve]{Stefan Samko\corauthref{cor}}
\corauth[cor]{Corresponding author.}
\ead{ssamko@ualg.pt}
\address[Algarve]{Universidade do Algarve, Campus de Gambelas,\\ Departamento de Matem\'atica, 8005-139 Faro, Portugal }
\ead[url]{w3.ualg.pt/\textasciitilde  ssamko}

\begin{abstract}
We prove that in variable exponent spaces $L^{p(\cdot)}(\Omega)$, where $p(\cdot)$ satisfies the
log-condition and  $\Omega$ is a bounded domain in $\mathbf R^n$  with the property that $\mathbf
R^n \backslash \overline{\Omega}$ has the cone property, the  validity of the Hardy type inequality
\[
\left \| \frac{1}{\delta(x)^\alpha}
 \int_\Omega \frac{\varphi(y)}{|x-y|^{n-\alpha}}\, dy \right \|_{p(\cdot)}  \leqq C \| \varphi
 \|_{p(\cdot)}, \quad 0<\al<\min\left(1,\frac{n}{p_+}\right)
\]
where $\delta(x)=\mathrm{dist}(x,\partial\Omega)$,  is equivalent to a certain property of the
domain $\Om$ expressed in terms of $\al$ and $\chi_\Om$.

\end{abstract}

\end{frontmatter}

\section{Introduction}

We consider  the Hardy inequality of the form
\begin{equation}\label{prin}
\left \| \frac{1}{\delta(x)^\alpha} \int_\Omega \frac{\varphi(y)}{|x-y|^{n-\alpha}}\, dy  \right
\|_{p(\cdot)}  \leqq C \| \varphi \|_{p(\cdot)}, \quad 0<\al<\min\left(1,\frac{n}{p_+}\right),
\end{equation}
within the frameworks of Lebesgue spaces with variable exponents $p(x),
p_+=\sup\limits_{x\in\Om}p(x)$, where $\delta (x)=\mathrm{dist} (x,\partial \Omega)$. We refer to
\cite{edmundsevans, edmundskokilmeskhi, kufmalpers} for Hardy type inequalities. The
multidimensional Hardy inequality of the form
\begin{equation}\label{necas}
\int_\Omega |u(x)|^p \delta (x)^{-p+a}dx \leqq C \int_\Omega |\nabla u(x)|^p \delta(x)^a dx,  \quad
u\in C^1_0(\Om),
\end{equation}
 appeared in \cite{nec} for bounded domains $\Omega \subset \mathbf R^n$ with Lipschitz boundary
and $1<p<\infty$ and $a>p-1$. This inequality was generalized by Kufner \cite[Theo. 8.4]{kuf} to
domains with H\"older boundary, and after that by Wannebo \cite{wan} to domains with generalized
H\"older condition. Haj\l asz \cite{haj} and Kinnunen and Martio \cite{kinmar} obtained a pointwise
inequality
\[
|u(x)|\leqq \delta(x)\mathcal M |\nabla u|(x),
\]
where  $\mathcal M$ is a kind of maximal function depending on the distance of $x$ to the boundary. This
 pointwise inequality combined with the
 knowledge of boundedness of Hardy-Littlewood maximal operator implies a ``local version near
 the boundary" of Hardy's inequality. This approach was used in the paper of
 Haj\l asz \cite{haj} in the case of classical Lebesgue spaces.

 Within the frameworks of variable exponent Lebesgue spaces, the Hardy inequality in one variable
  was first obtained in \cite {koksam04}, and later generalized in \cite{diesam}, where the
  necessary and sufficient conditions for the validity of the Hardy inequality on $(0,\infty)$ were
  obtained  under the assumption that  the log-condition  on $p(x)$ is satisfied only at the points
  $x=0$ and $x=\infty$, see also \cite{mascekmamogr, mascekogr}.

  For the multidimensional
  versions of Hardy inequality of form \eqref{prin} with  $\delta(x)^\alpha$ replaced by  $|x-x_0|^{\alpha},
  x_0\in\overline{\Omega}$, we refer to \cite{sam, samko}.
Harjulehto, H\"ast\"o and Koskenoja in \cite{harhaskos} obtained the  estimate
\[
 \left \| \frac{u(x)}{\delta (x)^{1-a}}\right\|_{p(\cdot)}\leqq C \| \nabla u(x) \delta(x)^a \|_{p(\cdot)},
 \quad u \in W^{1,p(\cdot)}_0(\Omega)
\]
making use of the approach of \cite{haj}, under the assumption that  $a$ is sufficiently small,
$0\leqq a<a_0$.

\vspace{3mm}Basing on some ideas and results of fractional calculus, in Theorem \ref{principal} we
show that the problem of the validity of inequality \eqref{prin} is equivalent to a certain
property of $\Om$ expressed in terms of $\al$ and $\chi_\Om$, see Definition \ref{bilateral} and
Theorem \ref{principal}. We did not find mentioning such an equivalence in the literature even in
the case of constant $p$.

\vspace{2mm} Note that the continuing interest to
 the
variable exponent  Lebesgue spaces  $L^{p(\cdot)}$ observed last
years was caused  by possible  applications (elasticity theory,
fluid mechanics, differential equations,  see for example
 \cite{525}).
We refer to papers \cite{kovrak91,618a} for basics on the
Lebesgue spaces with variable exponents and to the surveys
\cite{106b, 316b, 580bd} on harmonic analysis in
such spaces.
 One of the breakthrough results obtained for variable
$p(x)$ was the statement on the boundedness
 of the Hardy-Littlewood maximal operator in the generalized Lebesgue space $L^{p(\cdot)}$ under
certain conditions on $p(x)$, see \cite{106} and the further development in the above survey papers. The
importance of the boundedness of the maximal operator is known in particular due to the fact that many
convolution operators occurred in applications may be dominated by the maximal operator, which  is also used in
this paper.

Note also that the study of pointwise multipliers in the spaces of Riesz potentials is in fact  an
open question in case of variable $p(x)$.  Meanwhile, the topic of pointwise multipliers (in
particular, in the case of characteristic functions $\chi_\Om$) in spaces of differentiable
functions, is of importance in the theory of partial differential equations and other applications,
see for instance \cite{runstsickel}.

 The study of pointwise
multipliers of spaces of Riesz or Bessel potentials in the  case of constant $p$ may be found in
\cite{mazshap, mazshapR, str}, see also \cite{runstsickel} for the pointwise multipliers in the
case of more general spaces. We refer also, in the case of constant $p$ as well, to recent papers
\cite{662za, 662a} on the study of characteristic functions $\chi_\Om (x)$ as pointwise
multipliers.

\section{Preliminaries}

\subsection{On Lebesgue spaces with variable exponent}

The basics on variable Lebesgue spaces may be found in \cite{kovrak91,sam-diff}, but we recall here some
necessary definitions. Let $\Omega \subset \mathbf{R}^n$ be an open set. For a measurable function $p:\Omega
\rightarrow [1,\infty)$, we put
\[p^+=p^+(\Omega):=    \esss\limits_{x\in \Omega}p(x)   \;\; \ \
\textrm{and}\;\;\ \
 p^-=p^-(\Omega):=\essi\limits_{x\in \Omega}p(x)  .\]
In the sequel we use the  notation
\begin{equation}\label{bounds}
\mathcal{P}(\Omega) : = \{p\in L^\infty(\Omega): 1<p^-\leqq p(x) \leqq p^+<\infty\}.
\end{equation}

The generalised Lebesgue space $L^{p(\cdot)}(\Omega)$ with variable exponent is introduced as the set of all functions
$\varphi$ on $\Omega$ for which
\[\varrho_{p(\cdot)}(\varphi):=\int_{\Omega} |\varphi(x)|^{p(x)}dx < \infty.\]
Equipped with the norm
\[ \|\varphi \|_{L^{p(\cdot)}(\Omega)} := \inf \left\{ \lb > 0:
\varrho_{p(\cdot)}\left(\frac{\varphi}{\lb}\right)\leqq 1 \right\},\] this is a Banach space.
 The modular $\varrho_{p(\cdot)}(f)$ and the
norm $\|f \|_{p(\cdot)}$ are related to each other by
\begin{equation}\label{modularnorma}
\|f\|_{p(\cdot)}^\sigma \leqq I_p(f) \leqq \|f\|_{p(\cdot)}^\theta
\end{equation}

where $ \sigma=
\left\{\begin{array}{ll}\essi\limits_{x \in \Omega} p(x), & \|f\|_{p(\cdot)} \geqq 1; \\
\esss\limits_{x \in \Omega} p(x), & \|f\|_{p(\cdot)} \leqq 1\end{array}\right.
 $   and
$
 \theta=\left\{\begin{array}{ll}
\essi\limits_{x \in \Omega} p(x), & \|f\|_{p(\cdot)} \leqq 1; \\
\esss\limits_{x \in \Omega} p(x), & \|f\|_{p(\cdot)} \geqq 1.\end{array}\right. $

By $w$-$Lip\,(\Omega)$ we denote the class of all exponents $p\in L^\infty(\Omega)$ satisfying the
(local) logarithmic condition
\begin{equation}\label{llc}
    |p(x)-p(y)|\leqq \frac{C}{-\ln|x-y|}, \;\;\;\; |x-y|\leqq \frac{1}{2}, \ x,y \in \Omega.
\end{equation}

 By $p^\prime(\cdot)$ we denote the conjugate exponent, given by
$\displaystyle \frac{1}{p(x)}+\frac{1}{p^\prime(x)} \equiv 1.$

\subsection{Hardy-Littlewood maximal operator}

As usual, the Hardy-Littlewood maximal operator of a function $\vi$  on $\Om \subseteq \mathbf{R}^n$
is defined as
\begin{equation}
\mathcal M \varphi (x)=\sup_{r>0} \frac{1}{|\tilde{B}(x,r)|}\int_{\tilde{B}(x,r)}|\varphi(y)|dy,
\quad \tilde{B}(x,r)= B(x,r)\cap \Om.
\end{equation}

We use the notation
\begin{equation}\label{classP}
\mathbb{P}(\Omega) : = \left\{p: 1< p_-\leqq p_+\leqq \infty, \ \|\mathcal Mf\|_{L^{p(\cdot)}(\Omega)}\leqq
C\|f\|_{L^{p(\cdot)}(\Omega)}\right\}.
\end{equation}

\begin{prop}\emph{(\cite[Theo. 3.5]{106})}\label{harbou}
If $\Omega$ is bounded, $p\in \mathcal P(\Omega)\cap w$-$Lip(\Omega)$, then  $p\in
\mathbb{P}(\Om)$.
\end{prop}

\subsection{Potential and hypersingular integral operators}

\begin{defn}
For a function $\vi$  on $\mathbf{R}^n$, the Riesz potential operator $I^\alpha$ is defined by
\begin{equation}\label{riesz}
 I^\alpha \varphi (x)= \frac{1}{\gamma_n(\alpha)} \int\limits_{\mathbf R^n }
 \frac{\varphi(y)\, dy}{|x-y|^{n-\alpha}}
 = \varphi \ast k_\alpha(x),
 \end{equation}
where  the normalizing constant factor has the form $ \gamma_n(\alpha)=\frac{2^\alpha
\pi^{n\slash2} \Gamma\left( \frac{\alpha}{2}\right)}{\Gamma\left(\frac{n-\alpha}{2} \right)} $.
 The kernel
 $
k_\alpha(x)=\frac{|x|^{\alpha-n}}{\gamma_n(\alpha)}$ is referred to as the Riesz kernel.
\end{defn}

\begin{defn}
The space $I^\al\left(L^{p(\cdot)}\right)=I^\al\left(L^{p(\cdot)}(\mathbf{R}^n)\right), 0<\al<\frac{n}{p_+}$, called the
space of Riesz potentials,  is the space of functions $f$ representable as $f=I^\al\vi$ with
$\vi\in L^{p(\cdot)}, $ equipped with the norm
$\|f\|_{I^\al(L^{p(\cdot)})}=\|\vi\|_{L^{p(\cdot)}}$.
\end{defn}

\begin{defn}
The hypersingular integral operator $\mathbb D^\alpha$ of order $\al$, known also as the Riesz derivative,  is
defined by
\begin{equation}\label{lime}
\mathbb D^\alpha f = \lim_{\varepsilon \to 0} \mathbb D^\alpha_{\varepsilon} f =\lim_{\varepsilon
\to 0}  \frac{1}{d_{n,\ell}(\alpha)}\int\limits_{|y|>\varepsilon} \frac{\left(\Delta^\ell_y
f\right)(x)}{|y|^{n+\alpha}}dy ,
\end{equation}
where  $\alpha >0$ and $\ell>\al$ (see \cite[p.60]{samkobook},  for the value of the normalizing
constant $d_{n,\ell}(\al)$).
\end{defn}

It is known that given $\al$, one may choose an arbitrary   order $\ell>\al$ of the finite difference;   the
hypersingular integral does not depend on $\ell$ under this choice,  see  \cite[Ch. 3]{samkobook}.

In \cite{alm}, the following statement was proved.

\begin{prop}\label{hypinv}
Let $p\in \mathcal P(\mathbf R^n)\cap \mathbb P(\mathbf R^n)$ and $0<\al< \frac{n}{p_+}$. Then
\[ \mathbb D^\alpha I^\alpha \varphi=\varphi,\qquad \varphi \in L^{p(\cdot)}(\mathbf R^n) \]
where the hypersingular operator $\mathbb D^\alpha$ is taken in the sense of convergence of
$L^{p(\cdot)}\!$-norm.
\end{prop}

The characterization of the space $I^\al(L^{p(\cdot)}(\mathbf{R}^n))$ is given by the following
proposition.

\begin{prop} \emph{(\cite[Theo. 3.2]{almsam})} \label{rieszrange1}
Let $0< \alpha <n, p\in \mathcal P(\mathbf R^n)\cap \mathbb P(\mathbf R^n)$, $p^+<\frac{n}{\al}$ and let $f$ be a
locally integrable function.  Then $f\in I^\alpha \left( L^{p(\cdot)} \right)$ if and only if $f\in L^{q(\cdot)}$ with
$\frac{1}{q(\cdot)} = \frac{1}{p(\cdot)} - \frac{\alpha}{n}$, and there exists the Riesz derivative
$\mathbb{D}^\alpha f$ in the sense of convergence in $L^{p(\cdot)}$.
\end{prop}

\begin{rem}
Theorem 3.2 in \cite{almsam} was stated under the assumption that $p(x)$ satisfies the local
log-condition and the decay condition at infinity.  The analysis of the proof of Theorem 3.2 shows
that it is valid under the general assumption $p\in \mathcal P(\mathbf R^n)\cap \mathbb P(\mathbf
R^n)$ (if one takes into account that $p\in \mathcal{P}\cap\mathbb{P}(\mathbf{R}^n) \Leftrightarrow
p^\prime\in \mathcal{P}\cap\mathbb{P}(\mathbf{R}^n)$, see \cite[ Theo. 8.1]{106z}).
\end{rem}
 By Propositions \ref{hypinv} and \ref{rieszrange1},
for the norm $\|f\|_{I^\al(L^{p(\cdot)})}=\|\vi\|_{L^{p(\cdot)}}$ in the space of Riesz potentials
$I^\al\left(L^{p(\cdot)}(\mathbf{R}^n)\right)$ we have the following equivalence

\begin{equation} \label{equiv}
c_1\left(\|f\|_{L^{q(\cdot)}}+\|\mathbb{D}^\al f\|_{L^{p(\cdot)}}\right) \leqq
\|f\|_{I^\al\left(L^{p(\cdot)}\right)} \leqq c_2\left(\|f\|_{L^{q(\cdot)}}+\|\mathbb{D}^\al
f\|_{L^{p(\cdot)}}\right),
\end{equation}
where $\frac{1}{q(\cdot)} = \frac{1}{p(\cdot)} - \frac{\alpha}{n}$ and $c_1>0,c_2>0$ do not depend
on $f$.

\subsection{$(\al,p(\cdot))$-property of a domain  $\Om$}
\begin{defn}
A measurable function $g(x)$ is called a pointwise multiplier in the space
$I^\al\left(L^{p(\cdot)}(\mathbf{R}^n)\right)$, if $\|g I^\al\vi \|_{I^\al\left(L^{p(\cdot)}\right)}\leqq
C\|\vi\|_{L^{p(\cdot)}}$.
\end{defn}

By equivalence \eqref{equiv},  in the case $1<p_+<\frac{n}{\al}$ the characteristic function
$\chi_\Om(x)$  is a pointwise multiplier in $I^\al\left(L^{p(\cdot)}(\mathbf{R}^n)\right)$ if and only if
\begin{equation} \label{chimultiplier}
\|\mathbb{D}^\al (\chi_\Om I^\al \vi)\|_{L^{p(\cdot)}(\mathbf{R}^n)}\leqq
C\|\vi\|_{L^{p(\cdot)}(\mathbf{R}^n)} \quad \textrm{for all} \quad \vi \in
L^{p(\cdot)}(\mathbf{R}^n).
\end{equation}

 We introduce now the following notion related to the
property of the characteristic function $\chi_\Om$ to be a pointwise multiplier, but weaker than
that property. Let $\mathcal{E}_\Om f(x)=\widetilde{f}(x)=\left\{\begin{array}{ll} f(x),&  x\in\Om\\
0, & x\notin \Om
\end{array}\right.$ be the zero extension of a function $f$ defined on $\Om$.

\begin{defn}\label{bilateral}
We say that the domain $\Om$ has  the $(\al, p(\cdot))$-property, if  the function $\chi_\Om(x)$
has the following multiplier property
\begin{equation} \label{alphaproperty}
\|\mathbb{D}^\al(\chi_\Om I^\al \mathcal{E}_\Om \vi)\|_{L^{p(\cdot)}(\Om)}\leqq
C\|\vi\|_{L^{p(\cdot)}(\Om)} \quad \textrm{for all} \quad \vi \in L^{p(\cdot)}(\Om).
\end{equation}
\end{defn}

\begin{defn}
Let $p \in \mathcal P(\Om) \cap w$-$Lip (\Om)$. For brevity we call an extension $p^\ast(x)$ of
$p(x)$ to $\mathbf{R}^n$ regular, if $p^\ast \in \mathcal P(\mathbf R^n)\cap \mathbb P(\mathbf
R^n)$,  and $p^+(\mathbf{R}^n)=p^+(\Om)$. Such an extension is always possible, see \cite[Th.
4.2]{105a}; \cite[Lemma 2.2]{rafsam2}.
\end{defn}
\begin{lem}\label{implica}
 Let   $p \in \mathcal P(\Om)
\cap w$-$Lip (\Om)$. If $\chi_\Om$ is a pointwise multiplier in the space
$I^\al\left(L^{p^\ast(\cdot)}(\mathbf{R}^n)\right)$ under any regular extension $p^\ast(x)$ of $p(x)$ to
$\mathbf{R}^n$, then  the domain $\Om$ has the $(\al,p(\cdot))$-property.
\end{lem}
\begin{pf} We have to check condition \eqref{alphaproperty}, given that
$\|\chi_\Om f\|_{I^\al\left(L^{p^\ast(\cdot)}(\mathbf{R}^n)\right)}\leqq C
\|f\|_{I^\al\left(L^{p^\ast(\cdot)}(\mathbf{R}^n)\right)} $ under some regular extension of the exponent. We
have
$$\|\mathbb{D}^\al(\chi_\Om I^\al \mathcal{E}_\Om \vi)\|_{L^{p(\cdot)}(\Om)}\leqq
\|\mathbb{D}^\al(\chi_\Om I^\al \mathcal{E}_\Om \vi)\|_{L^{p^\ast(\cdot)}(\mathbf{R}^n)}.
$$
Since the extension $p^\ast(x)$ is regular, equivalence \eqref{equiv} is applicable so that
\begin{equation*}
\begin{split}
\|\mathbb{D}^\al(\chi_\Om I^\al \mathcal{E}_\Om \vi)\|_{L^{p(\cdot)}(\Om)} &\leqq C \|\chi_\Om
I^\al \mathcal{E}_\Om \vi\|_{I^\al\left(L^{p^\ast(\cdot)}(\mathbf{R}^n)\right)}\\
&\leqq C
 \|\mathcal{E}_\Om
\vi\|_{L^{p^\ast(\cdot)}(\mathbf{R}^n)}= C \| \vi\|_{L^{p(\cdot)}(\Om)},
\end{split}
\end{equation*}
which
completes the proof.
\end{pf}

\section{The main result}

\begin{thm}\label{principal}
Let $\Omega$ be a bounded domain in $\mathbf R^n$, $p \in \mathcal P(\Omega) \cap w$-$Lip (\Omega)$
and $0<\al<\min\left(1, \frac{n}{p_+}\right)$. If the domain $\Om$ has  the $(\al,
p(\cdot))$-property, then the Hardy inequality
\begin{equation}\label{princ}
\left \| \frac{1}{\delta(x)^\alpha} \int_\Omega \frac{\varphi(y)}{|x-y|^{n-\alpha}}\, dy \right
\|_{p(\cdot)}  \leqq C \| \varphi \|_{p(\cdot)}
\end{equation}
holds. If the exterior  $\mathbf R^n \backslash \overline{\Omega}$ has the cone property, then the
$(\al, p(\cdot))$-property is equivalent to the validity of the Hardy inequality \eqref{princ}.
\end{thm}

\section{Proof of Theorem \ref{principal}}

\subsection{The principal idea of the proof}

The proof of Theorem \ref{principal} is based on the observation that the weight
$\frac{1}{\delta(x)^\alpha}$ in fact is equivalent to the integral
$$a_\Om(x):=\int_{\mathbf{R}^n\backslash \Om}\frac{dy}{|x-y|^{n+\al}}, \quad x\in\Om.$$
Namely, the following statement is valid, see \cite[Prop. 3.1]{rafsam1}.

\begin{prop}\label{prop}
For an arbitrary domain $\Om$ there exists  a constant $c_1>0$ (not depending on $\Om,
c_1=\frac{1}{\al}|S^{n-1}|$) such that $a_\Om(x)\leqq \frac{c_1}{[\dl(x)]^\al}$. If the exterior
$\mathbf{R}^n\backslash \overline{\Om}$ has the cone property, then there exists a constant $c_2=
c_2(\Om)$ such that $\frac{1}{[\dl(x)]^\al}\leqq c_2 a_\Om(x)$.
\end{prop}

We will prove the following version of Theorem \ref{principal}.

\begin{thm}\label{principalc}
Let $\Omega$ be a bounded domain in $\mathbf R^n$, $p \in \mathcal P(\Omega) \cap w$-$Lip (\Omega)$
and $0<\al<\min\left(1, \frac{n}{p_+}\right)$.  Then the Hardy type inequality
\begin{equation}\label{princc}
\left \| a_\Om(x) \int_\Omega \frac{\varphi(y)}{|x-y|^{n-\alpha}} dy \right \|_{p(\cdot)}  \leqq C
\| \varphi \|_{p(\cdot)}
\end{equation}
holds if and only if the domain $\Om$ has the $(\al,p(\cdot))$-property.
\end{thm}

Theorem \ref{principal} will immediately follow from  Theorem \ref{principalc} in view of
Proposition \ref{prop}.
\subsection{On a hypersingular integral related to $\Om$.}


As in \cite{rafsam1}, we define the hypersingular integral (fractional derivative) of order
$0<\alpha<1$, related to the domain $\Om$, as the hypersingular integral over $\mathbf{R}^n$ of the
extension $\mathcal{E}_\Om f$:
$$\mathbb{D}_\Om f(x): =r_\Om \mathbb{D}^\al \mathcal{E}_\Om f(x)=\frac{1}{d_{n,1}(\al)}
\intl_{\mathbf{R}^n}\frac{f(x)-\widetilde{f}(y)}{|x-y|^{n+\al}}dy, \quad x\in \Om ,$$
where $r_\Om
$ stands for the restriction on $\Om$.  Splitting the integration in the last integral to that over
$\Om$ and $\mathbf{R}^n\backslash\ol{\Om}$, we can  easily see that
\begin{equation}\label{marchaud}
 a_\Omega(x)f(x) = d_{n,1}(\al)\mathbb{D}^\alpha \mathcal{E}_\Om f(x)- \int\limits_\Omega \frac{f(x)-f(y)}{|x-y|^{n+\alpha}} dy , \quad  x\in
\Omega.
\end{equation}

The proof of Theorem \ref{principalc} will be based on representation \eqref{marchaud} and certain
known facts from the theory of hypersingular integrals \cite{samkobook}.

\subsection{Auxiliary functions}

Although we will use the auxiliary functions defined below only in the case $\ell=1$, we give them
for an arbitrary integer $\ell$ as they are presented in \cite{samkobook}. By $\left(\Delta^\ell_h
f\right)(x):=\sum_{k=0}^\ell (-1)^k \binom{\ell}{k}f(x-kh)$ we denote the non-centered difference
 of a function $f$ defined on $\mathbf{R}^n$.
We need the non-centered difference
\begin{equation}
\Delta_{\ell,\alpha}(x,h):=\left(\Delta^\ell_h k_\alpha \right)(x)
\end{equation}
of the Riesz kernel  $k_\alpha(x)$ and single out the case of the step $h=e_1=(1,0,\ldots,0)$:
\begin{equation}
k_{\ell,\alpha}(x):=\Delta_{\ell,\alpha}(x,e_1) = \frac{1}{\gamma_n(\alpha)}\sum_{k=0}^\ell (-1)^k
\binom{\ell}{k}|x-ke_1|^{\alpha-n}.
\end{equation}

We  will also use  the function

\begin{equation}
\mathcal
K_{\ell,\alpha}(|x|)=\frac{1}{d_{n,\ell}(\alpha)|x|^n}\int\limits_{|y|<|x|}k_{\ell,\alpha}(y)dy.
\end{equation}

The following lemmata can be found in \cite[\S3.2.1]{samkobook}
\begin{lem}
The function $\Delta_{\ell,\alpha}(x,h)$, may be represented via its particular case
$k_{\ell,\alpha}(x)$ in terms of rotations:

\begin{equation}\label{rotequiv}
  \Delta_{\ell,\alpha}(x,h)=|h|^{\alpha-n} k_{\ell,\alpha} \left( \frac{|x|}{|h|^2}   \rota_x^{-1}h  \right)
\end{equation}

where $\rota_x\, \eta, \; \eta \in \mathbf R^n$ denotes any rotation in $\mathbf R^n$ which
transforms $\mathbf R^n$ onto itself so that $  \rota_x\; e_1=\frac{x}{|x|}. $
\end{lem}

\begin{lem}\label{3.13}
The function $k_{\ell,\alpha}(x)$ satisfies the condition
\begin{equation}\label{3.48}
|k_{\ell,\alpha}(x)| \leqq c(1+|x|)^{\alpha-n-\ell}\quad
\mathrm{when} \quad |x|\geqq \ell+1.
\end{equation}
\end{lem}

\begin{lem}\label{3.15}
Let $\ell>\Re \alpha>0$. Then
\begin{equation}
\int\limits_{\mathbf R^n} k_{\ell,\alpha}(y)dy=0.
\end{equation}
Moreover, in the case when $\ell$ is odd and the difference defining $k_{\ell,\alpha}(x)$ is non-centered,
\begin{equation}\label{3.55}
\int\limits_{\left|y-\frac{\ell}{2}e_1\right|<N}
k_{\ell,\alpha}(y)dy=0
\end{equation}
for any $N>0$.
\end{lem}

\begin{lem}
The function $\mathcal K_{\ell,\alpha}(|x|)$, $0<\alpha<1$ has the bound
\begin{equation}
| \mathcal K_{\ell,\alpha}(|x|)|\leqq C|x|^{\alpha-n}\quad \text{as }|x|\leqq 1.
\end{equation}

\end{lem}

\subsection{Proof  of Theorem \ref{principalc}}

Let $\vi\in L^{p(\cdot)}(\Om)$ and $\widetilde {\vi}=\mathcal{E}_\Om\vi (x)$. Substituting

$$f(y) : =I^\al \widetilde{\vi}= \frac{1}{\gm_n(\al)}
\intl_\Om\frac{\vi(t) }{|t-y|^{n-\al}}\,dt, \quad y\in\mathbf{R}^n$$
 into \eqref{marchaud}, we have
\begin{equation}\label{new}
a_\Omega(x)I^\alpha_\Omega \varphi (x)= \mathbb D^\alpha\chi_\Omega I^\alpha
\mathcal{E}_\Om\varphi(x)- \mathbb{A}\vi (x), \quad x\in \Om,
\end{equation}
where
$$\mathbb{A}\vi=\int\limits_\Omega \frac{I^\alpha \tilde{\varphi}(x)-I^\alpha \tilde{ \varphi}(y)}{|x-y|^{n+\alpha}}
dy = \lim\limits_{\ve\to 0}\mathbb{A}_\ve \vi(x)$$ and
$$ \mathbb{A}_\ve \vi(x)=
\int\limits_{ \substack{  y\in \Omega \\ |x-y|>\ve} }\frac{I^\alpha \tilde{\varphi}(x)-I^\alpha
\tilde{\varphi}(y)}{|x-y|^{n+\alpha}} dy.$$

The $(\al,p(\cdot))$-property of $\Om$, by the definition of this property and equivalence in \eqref{equiv},  is
nothing else but the boundedness in $L^{p(\cdot)}(\Om)$ of the operator $\mathbb D^\alpha\chi_\Omega I^\alpha
\mathcal{E}_\Om$. Thus, in the case of bonded domains $\Omega$,  the required equivalence of the Hardy inequality to the $(\al,p(\cdot))$-property will
follow from \eqref{new}, if the operator $\mathbb{A}$ is bounded.

\begin{lem}\label{hardom}
Let $0<\al<1$ and $\Om$ be a bounded domain.  The operators $\mathbb{A}_\ve $ are uniformly
dominated by the maximal operator:
\begin{equation}\label{eq2}
| \mathbb{A}_\ve \vi(x)|\leqq C \mathcal M\vi(x), \quad x\in \Om,
\end{equation}
for any $\vi\in L^1(\Om)$, where $C>0$ does not depend on $x$ and $\ve$. Consequently, the operator
$\mathbb{A}$ is bounded in the space $L^{p(\cdot)}(\Om)$ whenever $p\in \mathbb{P}(\Om)$.
\end{lem}
\begin{pf} We make use of the known representation
$$I^\alpha \tilde{\varphi}(x)-I^\alpha \tilde{\varphi}(x-y)=\int\limits_{\mathbf R^n}\Delta_{1,\al}(\xi,y) \tilde{\vi}(x-\xi) d\xi$$
 for the differences of the Riesz potential, see \cite[formula (3.64)]{samkobook},  and
get
\begin{eqnarray}\label{eq1}
\mathbb{A}_\ve \vi(x) &=&  \int\limits_{ \substack{ y\in\Omega_x \\ |y|>\ve}} \frac{dy}{|y|^{n+\alpha}} \int\limits_{\mathbf R^n} \tilde{\varphi}(x-\xi) \Delta_{1,\alpha}(\xi,y)d\xi  \nonumber\\
&=&\int\limits_{\mathbf R^n} \tilde{\varphi} (x-\xi)d\xi  \int\limits_{ \substack{ y\in\Omega_x \\ |y|>\ve}}
\frac{\Delta_{1,\alpha}(\xi,y)}{|y|^{n+\alpha}}dy
\end{eqnarray}
where $\Omega_x=\{y\in \mathbf R^n \colon x-y \in \Omega  \}$, {the interchange of the order of
integration being easily justified by Fubini's theorem whenever $\ve >0$. By \eqref{rotequiv} we
then have
\begin{eqnarray}\label{xi}
 \mathbb{A}_\ve \vi(x) &=& \int\limits_{\mathbf R^n} \tilde\varphi (x-\xi)d\xi
   \!\!\!\! \int\limits_{\substack{ y\in \Om_x \\ |y|>\varepsilon}}\frac{k_{1,\alpha}
   \left(  \frac{|\xi|}{|y|^2 }  \rota_\xi^{-1} \;  y \right)}{|y|^{2n}}dy  \nonumber\\
&=& \int\limits_{\mathbf R^n} \frac{ \tilde{\varphi} (x-\xi)}{|\xi|^n}d\xi \!\!\!\!
\int\limits_{\substack{ z\in\Om(x,\xi) \\ |z|<\frac{|\xi|}{\varepsilon}}}
k_{1,\alpha}(z)dz\nonumber\\ &=& \int\limits_{\mathbf R^n} \frac{ \tilde{\varphi}
(x-\ve\xi)}{|\xi|^n}d\xi \!\!\!\! \int\limits_{\substack{z\in\Om(x,\ve\xi) \\ |z|<|\xi|}}
k_{1,\alpha}(z)dz = \int\limits_{\mathbf R^n}  \tilde{\varphi} (x-\ve\xi)V_\varepsilon (x,\xi)d\xi,
\end{eqnarray}
where
\[\Om(x,\xi)=
\left\{ z\in \mathbf R^n \colon |\xi|\rota_\xi \frac{z}{|z|^2} \in \Om_x \right\}  \]}
and we denoted
$$V_\epsilon (x,\xi)=\frac{1}{|\xi|^n}\int\limits_{\substack{z\in\Om(x,\ve\xi)
\\ |z|<|\xi|}}  k_{1,\alpha}(z)dz$$
for brevity.
We  split $\mathbb A_\ve \vi(x)$ in the following way
\begin{equation}\label{split}
\mathbb A_\ve \vi(x)= \left( \int_{|\xi|<2} +\int_{|\xi|>2} \right) \tilde{\varphi} (x-\ve\xi) V_\varepsilon
(x,\xi)d\xi = : J_{1,\ve}\vi(x)+J_{2,\ve}\vi(x).
\end{equation}
For $J_{1,\ve}\vi(x)$ we have
\begin{eqnarray}\label{split1}
\left |J_{1,\ve}\vi(x) \right | &\leqq& \int\limits_{|\xi|<2} |\tilde{\varphi} (x-\ve\xi)|\frac{d\xi}{|\xi|^n} \int\limits_{|z|<|\xi|} |k_{1,\alpha}(z)|dz \nonumber \\
&\leqq& C  \int\limits_{|\xi|<2} \frac{ |\tilde{\varphi} (x-\ve\xi)|}{|\xi|^{n-\alpha}}d\xi \nonumber\\
&=& C \; |\tilde \vi| \ast \psi_\ve(x)
\end{eqnarray}
where $
\psi(\xi)= \left\{\begin{array}{ll} |\xi|^{\alpha-n}, & |\xi|<2, \\
0, & |\xi| \geqq 2, \end{array}\right. $ and $\psi_\ve(x)=\ve^{-n}\psi(x/\ve)$.

When $|\xi|>2$, the key moment in the estimation is the usage of property \eqref{3.55} of the Riesz kernel:
\begin{eqnarray*}
V_\ve(x,\xi) & \substack{ \\ =} & \frac{1}{|\xi|^n}
 \left(\;  \int\limits_{B(0,|\xi|) \cap  \Omega(x,\ve\xi) } -
 \int\limits_{\left|z-\frac{e_1}{2}\right|<|\xi|-1}  \; \right)  k_{1,\alpha}(z)dz \\
&=&\frac{1}{|\xi|^n}\int\limits_{\Theta(x,\ve)}k_{1,\alpha}(z)dz
\end{eqnarray*}

where
\[
\Theta(x,\ve)=\left\{z:z\in B(0,|\xi|) \cap \Omega(x,\ve\xi)
\right\} \backslash \left\{z:\left|z-\frac{e_1}{2}\right|<|\xi|-1
\right\}.
\]

Since $\Theta(x,\ve)$ is embedded in the annulus
$|\xi|-\frac{3}{2}\leqq |z|\leqq |\xi|$, we have

\begin{eqnarray*}
|V_\ve(x,\xi)| &\leqq&
\frac{1}{|\xi|^n}\int\limits_{|\xi|-\frac{3}{2}\leqq |z|\leqq
|\xi|} |k_{1,\alpha}(z)|dz
\end{eqnarray*}
and by \eqref{3.48}

\begin{eqnarray}\label{3.50}
|V_\ve(x,\xi)| &\leqq& \frac{C}{|\xi|^n}\left| |\xi|^{\alpha-1}-
\left(|\xi|-\frac{3 }{2}\right)^{\alpha-1} \right| \leqq
\frac{C}{|\xi|^{n+2-\alpha}}.
\end{eqnarray}

The estimation of $J_{2,\ve}\vi(x)$ is then given by

\begin{equation}\label{split2}
|J_{2,\ve}\vi(x)| \substack{\eqref{3.50}\\ \leqq} C \; |\tilde \vi|
\ast \phi_\ve(x)
\end{equation}
 where
$
\phi(\xi)= \left\{\begin{array}{ll} 2^{\alpha-n-2}, & |\xi|<2, \\
|\xi|^{\alpha-n-2}, & |\xi| \geqq 2, \end{array}\right. $ and $\phi_\ve(x)=\ve^{-n}\phi(x/\ve)$.

Since the kernels $\psi$, $\phi$ are radially decreasing and  integrable, we can use the well known estimation of
convolutions with such  kernels via the maximal function, which yields
\begin{equation}\label{maximal}
J_{i,\ve}\vi(x)\leqq C\mathcal M(|\varphi|), i=1,2, \forall \ve>0
\end{equation}
and implies \eqref{eq2} after gathering \eqref{split}, \eqref{split1}, \eqref{split2} and \eqref{maximal}. This
completes the proof.
\end{pf}
\subsection{Corollaries}

As a corollary of Theorem \ref{principal} we obtain an estimate in classical $L^p(\Omega)$ spaces, but first  we
need the following definition.

\begin{defn}
Let $\Omega$ be an open set in $\mathbf{R}^n$. We say that
$\Omega$ satisfies the  \textit{Strichartz condition} if  there
exist a coordinate system in $\mathbf{R}^n$ and an integer $N>0$
such that almost every line parallel to the axes intersects $\Omega$
in at most $N$ components.
\end{defn}

\begin{lem}\emph{(\cite{nogrub,str};\cite[p. 244]{rubin})}. \label{pointwise}
The characteristic function $\chi_\Omega$ of a domain $\Omega$
satisfying the Strichartz condition is a pointwise multiplier in
the space $I^\alpha\left( L^p(\mathbf R^n) \right)$  when
$1<p<1/\alpha$.

\end{lem}

\begin{cor}
The Hardy inequality
\[
\left \| \frac{1}{\delta(x)^\alpha}\int\limits_{\Omega} \frac{\vi(y)}{|x-y|^{n-\alpha}}dy
\right\|_p \leqq C \|\vi\|_p, \quad 1<p<1/\alpha
\]
holds for any bounded open set  $\Omega\subset \mathbf R^n$ satisfying the  Strichartz condition.
\end{cor}

\begin{pf}
By Lemma \ref{implica} and Lemma \ref{pointwise} we have that
$\Omega$ has the $(\alpha,p(\cdot))$-property and then the results
follows from Theorem \ref{principal}.
\end{pf}







\end{document}